\magnification=\magstep0
\input amstex

\documentstyle{amsppt}
\NoBlackBoxes
\NoRunningHeads
\nologo
\topmatter
\title{A matrix equation for association schemes}
\endtitle
\author L. Chihara
and D. Stanton
\footnote{\hbox{This work was supported by NSF grant DMS-9400510.\hfill}}
\endauthor
\address Department of Mathematics, St. Olaf College,
Northfield, MN 55057.
\endaddress
\address School of Mathematics, University of Minnesota,
Minneapolis, MN 55455.
\endaddress
\abstract{Let ${\Cal X}$ be a self-dual P-polynomial association scheme. 
Then there are at most 12 diagonal matrices $T$ such that $(PT)^3=I$.
Moreover, all of the solutions for the classical 
infinite families of such schemes 
(including the Hamming scheme) are classified.}  
\endabstract
\endtopmatter

\subheading{1. Introduction} 

Let ${\Cal X}=(X,\{R_i\})$ be a self-dual P-polynomial association scheme. 
In \cite{2}, a spin model is constructed 
when ${\Cal X}=H(d,q)$, the Hamming scheme. It relies on a solution to
the matrix equation 
$$
(PT)^3=I
\tag 1.1
$$
where $T$ is a diagonal matrix. 
In this paper we prove that (1.1) has no other solutions for the Hamming 
scheme, other than the ones given in \cite{2}. We also prove that a general 
P-polynomial association scheme has at most 12 solutions to (1.1). We 
classify the solutions to (1.1) for the self-dual matrix schemes 
over a finite field.

\subheading{2. The general case} 

We assume that ${\Cal X}$ is an $N$-class self-dual 
$P$-polynomial association scheme.

In this section we derive some necessary conditions on the entries 
of $T$ for (1.1) to hold, and give a proof of the following theorem.

\proclaim{Theorem 1} If ${\Cal X}$ is an $N$-class self-dual 
P-polynomial association scheme, with $N\ge 2$, then there are at 
most 12 solutions $T$ to $(PT)^3=I$. 
\endproclaim 

In fact, we will show that there are at most 4 values for 
$T_1/T_0$ (given by solutions to (2.6) below), at most 3 values for 
$T_0$, and all other entries of $T$ are uniquely determined from $T_0$ 
and $T_1$. 

We shall use notation for association schemes found in \cite{4}, 
except that we let $v_i$ (not $k_i$) denote the 
$ith$ valency. Since ${\Cal X}$ is self-dual, $v_i=m_i$, the $ith$ 
multiplicity, and $P_j(i)=Q_j(i)$.

\demo{Proof of Theorem 1}
Let 
$T=diag(T_0,T_1,\dots,T_N)$. Since $(PT)^2=(PT)^{-1}$, by 
finding the $(i,j)$ entry of both sides we have
$$
{P_j(i)\over |X|T_i T_j}=\sum_{k=0}^N T_k P_k(i)P_j(k).
\tag2.1a
$$
Similarly, setting $(PT)^{-2}=PT$, we find
$$
|X|T_iT_jP_j(i)=\sum_{k=0}^N {P_k(i)P_j(k)\over T_k}.
\tag2.1b
$$

The intersection numbers $p_{i,j}^r$ of the scheme satisfy \cite{4 ,p. 53}
$$
P_i(k)P_j(k)=\sum_{r=0}^N p_{i,j}^r P_r(k).
\tag 2.2
$$
Since \cite{4, p. 63} 
$$
{P_k(i)\over v_k}={Q_i(k)\over m_i}={P_k(i)\over m_i},
$$
we can expand the right sides of (2.1a) and (2.1b) using (2.2). 
The results are
$$
{T_0v_iP_j(i)\over T_iT_j}=\sum_{r=0}^N {p_{i,j}^r v_r \over T_r}
\tag2.3a
$$
and 
$${v_iT_iT_jP_j(i)\over T_0}=\sum_{r=0}^N p_{i,j}^r v_rT_r
\tag2.3b
$$

Recall that $p_{i,j}^r=0$ if $|r-j|\ge 2$;
also, $c_j=p_{1,j-1}^{j}, a_j=p_{1,j}^j, b_j=p_{1,j+1}^j$, and
$v_j=\prod_{i=0}^{j-1}\frac{b_i}{c_{i+1}}$. 
Putting $j=1$ in (2.3a) and (2.3b) gives
$$
{v_iT_iT_1P_1(i)\over T_0T_0}=b_{i-1}v_{i-1}{T_{i-1}\over T_0}+a_iv_i{T_i\over T_0}+
c_{i+1}v_{i+1}{T_{i+1}\over T_0},
\tag2.4a
$$
and 
$$
v_i{T_0T_0P_1(i)\over T_iT_1}=b_{i-1}v_{i-1}{T_0\over
T_{i-1}}+a_iv_i{T_0\over T_{i}} +c_{i+1}v_{i+1}{T_0\over T_{i+1}}.
\tag2.4b
$$

First we show that $T_i$ is uniquely determined from $T_1/T_0=x$ and $T_0$. 
If $P_1(i)\ne 0$ for all i, then it is clear from (2.4a) \cite{5} 
that $T_i/T_0$ is an 
orthogonal polynomial in $x$ of degree $i$, which we will call $t_i(x)$. The 
last equation in (2.4a) ($i=N$) gives a polynomial equation for $x$ 
of degree $N+1$. If $P_1(i)=0$ for some $i$, that $i$ must be unique, 
because the eigenvalues are distinct. Then (2.4a) implies that 
$T_{i+1}/T_0$ is a polynomial in $x$ of degree $i$, and clearly that
$T_{j+1}/T_0$ is a polynomial in $x$ of degree $\le j$, for $j\ge i$. 
In any case $T_i/T_0$ is uniquely determined from $T_1/T_0=x$ and $T_0$.  

Similarly, (2.4b) shows that $T_0/T_i=t_i(1/x)$. Thus 
we must have
$$
t_i(x)t_i(1/x)=1, \text{ for } 1\le i\le N.
\tag2.5
$$
Clearly, (2.5) is satisfied for $i=1$ since $t_1(x)=x$. However the case 
$i=2$ is a quartic equation for $x$,
$$
P_1(1)x^4+a_1(P_1(1)-1)x^3-(P_1(1)^2+a_1^2-b_1^2+1)x^2+a_1(P_1(1)-1)x
+P_1(1)=0.
\tag 2.6
$$
(2.6) is a non-trivial equation, since if $P_1(1)=0=a_1$ and $b_1=1$, 
we would have $k=2$, thus a $4$-gon.
So $x$ has at most 4 possible values. $T_0$ has at 
most 3 values, because it must satisfy a cubic equation from $(PT)^3=I$.
\qed\enddemo

Theorem 1 is best possible, since even $n$-gons have 12 solutions 
(see Theorem 6).
The cases $N=1$ and $N=2$, always have solutions to (1.1), see \cite{3}.
So we can assume that $N>2$.

\subheading{3. Classical cases}

In this section we apply the necessary conditions of \S2 to 
the infinite families classical self-dual P-polynomial association 
schemes in \cite{4}.

\proclaim{Theorem 2} The only solutions to $(PT)^3=I$ in the Hamming scheme
$H(N,q)$ for $N>2$, $q\ge 2$ 
are $T_i=c x^i$, where $x$ satisfies $1-2x+qx+x^2=0$ 
and $c$ satisfies $c^3(q(1+(q-1)x))^N=1$. 
\endproclaim
\demo{Proof} Here we have (see \cite{4}) 
$c_i=i$, $a_i=i(q-2)$, $b_i=(N-i)(q-1)$, 
$P_1(i)=N(q-1)-qi$, and $v_i=\binom{N}{i}(q-1)^i$. 
From \S2, we know that $T_i/T_0=t_i(x)$
is a polynomial of degree $i$ in $x$. So $x$ must be a simultaneous 
zero of the numerator of $t_2(x)t_2(1/x)-1$,
$$
(1 - 2x + qx + x^2)(N + q - Nq - 2x + qx + *x^2 + qx^2 - Nqx^2)
$$
and the numerator of $t_3(x)t_3(1/x)-1$,
$$
(1 - 2x + qx + x^2)(\text{a certain quartic in  } x).
$$
Suppose that $x$ does not satisfy $1-2x+qx+x^2=0$. The resultant of these two 
polynomials in $x$ divided by $1-2x+qx+x^2$ is
$$
  4(N-1)^2(q-2)^2(q-1)^2(-2 - N + Nq)^2 (-N - q + Nq)^4.
\tag3.1
$$
The resultant must be zero for a common zero. Since $N>2$ and $q\ge2$, 
(3.1) forces $q=2$. If $q=2$, the numerator of (2.4b) with $i=1$ is $N(N-2)(1+x^2)^2$, 
which implies that $x=\sqrt{-1}$\thinspace
-already the allowed zero. 

It is not hard to show that $t_i(x)=x^i$ by induction using (2.4a). For the 
normalization constant $c$, compute the $(0,0)$-entry of $(PT)^3=I$ 
and use the generating functions for the Krawtchouk polynomials. 
\cite{5}.
\qed\enddemo

\proclaim{Theorem 3} There are no solutions to $(PT)^3=I$ in the association 
scheme of $M \times N$ bilinear forms over a finite field of order $q$ if 
$min\{M,N\}>2$. 
\endproclaim
\demo{Proof}
In this scheme (see \cite{4}) if $d=q^M$, $e=q^N$, we have
$b_i=q^{2i}(dq^{-i}-1)(eq^{-i}-1)/(q-1)$, 
$c_i=q^{i-1}(q^i-1)/(q-1)$,
$a_i=(1-e)(1-d)/(q-1)-b_i-c_i$,
and $P_1(i)= v_1(de+q^i-dq^i-eq^i)/(d-1)(e-1)q^i$. We use (2.4b) directly, 
rather that $t_i(x)t_i(1/x)=1$. 
Using the explicit formula for $t_2(x)$ and $t_3(x)$
as polynomials in $x$,
we let $i=1$ and $i=2$ in (2.4b). We obtain two polynomials in $x$,$d$,$e$, 
and $q$ which must be zero.
The greatest common divisor this time is $(d-1)(e-1)$, which is non-zero, 
so the resultant, as polynomials in $d$, must be zero. It is
$$
\align
R=&(q-1)^5 q^{12}(e -q)^6(e -q^2)^4(x-1)^4x^3(1-2x+ex+x^2)^6\\
&\times (-1+e-q+2x-2ex+e^2x+2qx-2eqx-x^2+ex^2-qx^2)^2.\\
&\times(-q-x- q^2x+q^4x-qx^2)=0\\
\endalign
$$

Since we are assuming $e,d>q^2$, $x$ must be a root of one of the factors
in $R=0$. This gives four possible polynomials for which $x$ is a zero. 
We shall find the 
remainder when these four polynomials in $x$ divide (2.4b) with $i=1$. 
We will see that this remainder cannot be zero, giving a contradiction.

For $x-1$ the remainder is
$$
de(-de - deq + 2dq^2 + 2eq^2 - 2q^3)=0.
\tag3.3
$$
Since $q$ is a prime power, and $e$ and $d$ are both positive powers of $q$,
we have a unique smallest power of $q$ in (3.3), so there is no solution.

For $1-2x+ex+x^2=0$, the remainder is
$$
e(d - q)^2(e - q^2)(1 - 2x + ex)=0.
\tag3.4
$$
The only solution is $x=-1/(e-2)$, which then implies that 
$1-2x+ex+x^2=1/(2-e)^2\ne 0$.

For $-q-x -q^2x+q^4x-qx^2=0$, the remainder is
$$
\aligned
&(-1-2q+dq + eq - deq - 2q^2 + dq^2 + eq^2 - deq^2 - q^3 + dq^3 + eq^3)\\
&\times(-de-q+dq+eq-2q^2 + dq^2 + eq^2 - dq^3 - eq^3 + 2q^4 + q^5 - q^6)\\
&\times(q + x + q^2x - q^4x)/q^2=0.\\
\endaligned
\tag3.5
$$
Reducing modulo $q$ shows the first factor is non-zero, 
similarly the $-q$ term in the second factor shows it is non-zero. 
Thus we must have $x=q/(q^4-q^2-1)$, which again contradicts 
$-q-x -q^2x+q^4x-qx^2=0$.

For $-1+e-q+2x-2ex+e^2x+2qx-2eqx-x^2+ex^2-qx^2=0$, the remainder is
$$
\aligned
&(-de + 2de^2 - de^3 - 3deq - 2e^2q + 3de^2q + 
       e^3q + 2dq^2 + 5eq^2 - 4deq^2 - 3e^2q^2 \\
&- 2q^3 + 2dq^3 + 3eq^3 - 2q^4)e(-d + q^2)/(-1 + e - q)^3\\
&\times(-1 + e - q + 2x - 2ex + e^2x + 2qx - 2eqx)=0.
\endaligned
\tag3.6
$$
Again the $-2q^3$ term shows the first factor is non-zero, so we must have
$x=(1-e+q)/(2-2e+e^2+2q-2eq)$, which implies that
$$-1+e-q+2x-2ex+e^2x+2qx-2eqx-x^2+ex^2-qx^2=
\frac{(e-q-1)^2}{(2-2e+e^2+2q-2eq)^2}.
\tag 3.7
$$
Since $e>q+1$, (3.7) cannot be zero.

So in all cases, the resultant is non-zero, and the two polynomials in $d$ 
cannot have a common zero.
\qed\enddemo

A nearly identical proof gives non-existence for the 
alternating and Hermitian forms. Again we need at least 3 classes.
 
\proclaim{Theorem 4} There are no solutions to $(PT)^3=I$ in the association 
schemes of $N \times N$ alternating forms over a finite field of order 
$q$ if $N>5$. 
\endproclaim

\proclaim{Theorem 5} There are no solutions to $(PT)^3=I$ in the association 
schemes of $N \times N$ Hermitian forms over a finite field of order 
$q^2$ if $N>2$. 
\endproclaim

For completeness we state the result for the $n$-gons, see \cite{1}. 
We let $i=\sqrt{-1}$.

\proclaim{Theorem 6} The only solutions to $(PT)^3=I$ in the association 
schemes of regular $n$-gons for $n\ge6$ are 
$$
T_j=c(-1)^j e^{\pm\pi i j^2/n},
$$ 
where
$$
c^3 n^{3/2}(-1)^m =
\cases e^{\pm\pi i/4}\text{ if } n=4m\\
 e^{\mp\pi i/4}\text{ if } n=4m+2\\
1 \text{ if } n=4m+1\\
\mp i\text{ if } n=4m+3,
\endcases
$$
and for $n$ even 
$$
T_j=c e^{\pm\pi i j^2/n},
$$
where
$$
c^3 n^{3/2}=e^{\pm\pi i/4}.
$$
\endproclaim
\demo{Sketch of proof} The equation $t_2(x)t_2(1/x)=1$ has four solutions, 
$x=\pm e^{\pm\pi i/n}$. It is easy to see that 
$T_j=c(-1)^j e^{\pm\pi i j^2/n}$ for $x=-e^{\pm\pi i/n}$, and 
$T_j=c e^{\pm\pi i j^2/n}$ for $x=e^{\pm\pi i/n}$. However for $n$ odd, 
the latter solution does not satisfy (2.4a) for $i=(n-1)/2$. 
A Gauss sum calculation verifies the solutions given, and 
finds the constant $c$.
\qed\enddemo

\enddocument
\end